\documentclass{amsart}

\usepackage{xcolor}
\usepackage{graphicx}
\numberwithin{equation}{section}
\usepackage{cite}
\usepackage{amsmath,amssymb,mathtools}
\usepackage[hypertexnames=false]{hyperref}
\usepackage[initials,nobysame]{amsrefs}
\hypersetup{
	colorlinks=true,
	allcolors=red!70!black,
	pdftitle={Topping's diameter conjecture for immersed surfaces},
	pdfauthor={Tatsuya Miura},
	pdfsubject={Topping's diameter conjecture and sharp mean-width estimates for immersed surfaces},
	pdfkeywords={Mean curvature, diameter, mean width, immersed surfaces, minimal surfaces},
}
\newtheorem{theorem}{Theorem}[section]
\newtheorem{proposition}[theorem]{Proposition}
\newtheorem{lemma}[theorem]{Lemma}
\newtheorem{corollary}[theorem]{Corollary}
\theoremstyle{remark}

\theoremstyle{definition}
\newtheorem*{acknowledgements}{Acknowledgements}
\newtheorem*{AI}{AI usage disclosure}

\usepackage[capitalize]{cleveref}
\crefname{theorem}{theorem}{theorems}
\Crefname{theorem}{Theorem}{Theorems}
\crefname{proposition}{proposition}{propositions}
\Crefname{proposition}{Proposition}{Propositions}
\crefname{lemma}{lemma}{lemmas}
\Crefname{lemma}{Lemma}{Lemmas}
\crefname{corollary}{corollary}{corollaries}
\Crefname{corollary}{Corollary}{Corollaries}
\crefname{remark}{remark}{remarks}
\Crefname{remark}{Remark}{Remarks}

\newcommand{\R}{\mathbf R}
\newcommand{\Sph}{\mathbf S^2}
\newcommand{\HH}{\mathcal H}
\DeclareMathOperator{\dist}{dist}
\DeclareMathOperator{\diam}{diam}

\DeclareMathOperator{\tr}{tr}

\title[Topping's diameter conjecture for immersed surfaces]{Topping's diameter conjecture for immersed surfaces}
\author[T.~Miura]{Tatsuya Miura}
\address[T.~Miura]{Department of Mathematics, Graduate School of Science, Kyoto University, Kitashirakawa Oiwake-cho, Sakyo-ku, Kyoto 606-8502, Japan}
\email{tatsuya.miura@math.kyoto-u.ac.jp}
\date{\today}
\subjclass[2020]{53A05 (primary), 53A10 (secondary)}
\keywords{Mean curvature, diameter, mean width, immersed surfaces, minimal surfaces}

\begin{document}
\begin{abstract}
We resolve Topping's diameter conjecture for closed immersed surfaces in Euclidean three-space, building on his recent result for embedded surfaces. More generally, we establish sharp diameter and mean-width estimates for compact connected immersed surfaces, possibly with boundary. As a consequence, we confirm the codimension-one case of the sharp length--diameter conjecture for immersed minimal surfaces.
\end{abstract}

\maketitle


\section{Introduction}

In 1998, Topping \cite{Topping1998}*{(17)} conjectured that every closed connected surface immersed in $\R^3$ satisfies the sharp diameter bound
\[
	\diam(f(\Sigma))<\frac1\pi\int_\Sigma|\mathbf H_f|\,d\mu_f,
\]
where $f:\Sigma\to\R^3$ is the immersion, $d\mu_f$ is its induced area measure, and $\mathbf H_f$ is the mean curvature vector normalized so that the unit sphere has $|\mathbf H_f|=1$.
Here and below, diameter means extrinsic diameter $\diam(A)\coloneqq \sup_{x,y\in A}|x-y|$ for $A\subset\R^3$.

Building on Topping's resolution of the embedded case \cite{Topping2026}, we prove the conjecture for the general class of closed immersed surfaces and obtain its optimal extension to surfaces with boundary.

Throughout this paper, surfaces, their boundaries, and immersions are smooth.
We write $\operatorname{Length}(\partial\Sigma)\coloneqq \int_{\partial\Sigma}d\mu_{f|_{\partial\Sigma}}$ for boundary length with respect to the measure induced by $f|_{\partial\Sigma}$, with value zero when $\partial\Sigma=\varnothing$.
Area integrals and boundary lengths count multiplicity.
For a nonempty compact set $A\subset\R^3$, let
\[
W_A(e)\coloneqq \max_{x\in A}\langle x,e\rangle-\min_{x\in A}\langle x,e\rangle
\]
denote the width in direction $e\in\Sph$, where $\Sph\subset\R^3$ is the unit sphere.
Let $\HH^2$ denote the two-dimensional Hausdorff measure.

Our main result is a sharp mean-width estimate for immersed surfaces, from which the diameter bound follows.

\begin{theorem}\label{thm:main}
	Let $f:\Sigma\to\R^3$ be an immersion of a compact connected surface $\Sigma$, possibly with boundary $\partial\Sigma$.
	Then
	\begin{equation}\label{eq:main}
		\begin{aligned}
			\int_\Sigma|\mathbf H_f|\,d\mu_f+\frac\pi2\operatorname{Length}(\partial\Sigma)
			\geq\frac12\int_{\Sph}W_{f(\Sigma)}(e)\,d\HH^2(e)>\pi\diam(f(\Sigma)).
		\end{aligned}
	\end{equation}
	In particular,
	\[
		\diam(f(\Sigma))
		<\frac1\pi\int_\Sigma|\mathbf H_f|\,d\mu_f+\frac12\operatorname{Length}(\partial\Sigma).
	\]
\end{theorem}

The coefficients $1/\pi$ and $1/2$ in the diameter estimate are sharp already for embedded surfaces, as seen from long capped cylinders and long thin smooth planar domains, respectively.
The mean-width inequality is also optimal, with equality attained by closed convex surfaces and compact convex planar domains.

For closed convex surfaces, the diameter estimate is a classical consequence of Minkowski's integral formula \cite{Minkowski1903}.
For closed orientable constant-mean-curvature surfaces, Topping \cite{Topping1997}*{Theorem~5} proved a stronger diameter bound than the conjectured one.
The conjectured sharp bound was established for simply connected surfaces of revolution in \cites{Helmers2015,Miura2021}.
The author's work \cite{Miura2021} also gives a quantitative remainder and a description of near-equality sequences.

Related estimates for general closed surfaces include Simon's bound in terms of area and Willmore energy \cite{Simon1993}.
Topping \cite{Topping2008}*{Theorem~1.1} bounded the intrinsic diameter of closed immersed submanifolds of arbitrary dimension and codimension in terms of mean curvature, with a constant depending only on the dimension.

For surfaces with boundary, Paeng \cite{Paeng2014} obtained intrinsic diameter estimates under a geodesic convexity assumption.
Without this assumption, Menne--Scharrer \cite{MenneScharrer2026}*{Corollary~2} established estimates in Euclidean spaces through a general varifold theory, while the author \cite{Miura2022} developed a geometric method aimed at obtaining sharp constants for surfaces.
Related diameter estimates in curved ambient spaces, including results for submanifolds with boundary, were obtained in \cites{WuZheng2011,Wu2023,FlaimScharrer2024,LiuTaiWu2025,ChambersMarxKuo2026}.
For applications of diameter estimates in other geometric contexts, we refer to \cite{Topping2026} and the references therein.

Recently, Topping \cite{Topping2026}*{Theorem~1.4} proved that every connected closed embedded surface $M\subset\R^3$ satisfies
\begin{equation}\label{eq:topping}
	\frac12\int_{\Sph}W_M(e)\,d\HH^2(e)\leq\int_M H_+\,d\HH^2.
\end{equation}
Here $H\coloneqq\frac12\operatorname{div}_M\nu$ for the outward unit normal $\nu$, and $H_+\coloneqq\max\{H,0\}$.
His proof uses Spadaro's mean-convex hull \cite{Spadaro2020} to fill holes in planar sections, and then applies a Crofton formula for total mean curvature.
In particular, this proves the sharp diameter conjecture for embedded closed surfaces.
For compact connected embedded surfaces with boundary, Topping \cite{Topping2026}*{Theorem~1.8} also proved the mean-width and diameter estimate \eqref{eq:main}, with the boundary coefficient $\pi/2$.
Our result extends the latter estimate from embeddings to immersions.

For minimal surfaces, \Cref{thm:main} immediately gives the following sharp length--diameter estimate.

\begin{corollary}\label[corollary]{cor:minimal}
	Let $f:\Sigma\to\R^3$ be a minimal immersion of a compact connected surface with boundary.
    Then
	\[
		\operatorname{Length}(\partial\Sigma)>2\diam(f(\Sigma))
		=2\diam(f(\partial\Sigma)).
	\]
\end{corollary}

\begin{proof}
	The strict inequality follows from \Cref{thm:main} because $\mathbf H_f=0$.
	The diameter equality follows since $f(\Sigma)$ lies in the convex hull of $f(\partial\Sigma)$ by the maximum principle.
\end{proof}

The optimal length--diameter constant $2$ for immersed minimal surfaces was conjectured by the author in
\cite{Miura2022}*{(1.2)} for all ambient dimensions $n\geq3$.
The above corollary confirms this conjecture in codimension one, extending the corresponding result for embedded minimal surfaces by Topping
\cite{Topping2026}*{Corollary~1.9}.
In particular, this result gives a sharp nonexistence criterion for the Plateau--Douglas problem, as discussed in \cites{Miura2022,Topping2026}.

\begin{figure}[htbp]
	\centering
	\begin{minipage}[t]{.48\linewidth}
		\centering
		\includegraphics[width=\linewidth]{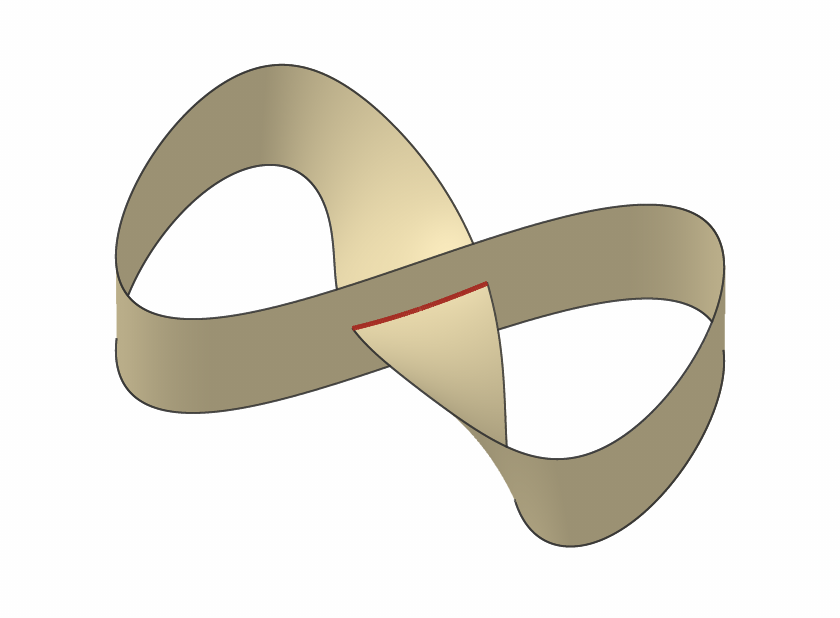}
		\par\smallskip
		\normalfont (a)
	\end{minipage}\hfill
	\begin{minipage}[t]{.48\linewidth}
		\centering
		\includegraphics[width=\linewidth]{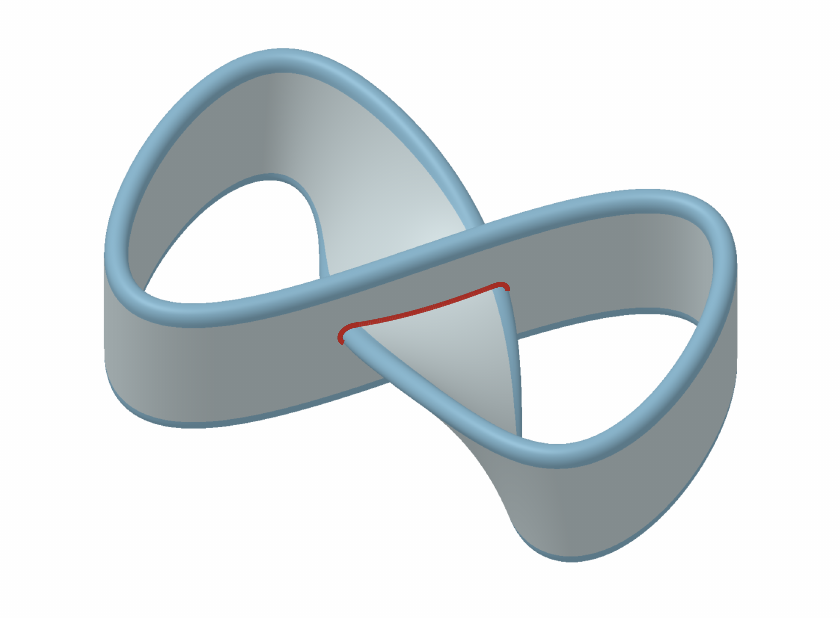}
		\par\smallskip
		\normalfont (b)
	\end{minipage}
	\caption{(a) An immersed M\"obius band, which is nonorientable.
	(b) The orientable closed surface obtained by doubling.
	The black curve in (a) marks the original boundary; blue rounded collars in (b) join the gray surfaces on either side of the original surface smoothly.
	Red curves indicate the visible self-intersections.}
	\label{fig:doubling-sequence}
\end{figure}

Our proof reduces the general problem to Topping's estimate \eqref{eq:topping} for closed embedded surfaces.
We use the doubling trick from the author's previous work \cite{Miura2022}, also used by Topping \cite{Topping2026}, to construct a closed comparison surface.
In this construction, the two local normal perturbations of $f$ approximate the original surface away from its boundary, while narrow rounded collars join them along the boundary, as illustrated in \Cref{fig:doubling-sequence}.
These collars behave like curved half-cylinders: their cross-sectional profiles turn through angle $\pi$, yielding the limiting boundary contribution $\frac\pi2\operatorname{Length}(\partial\Sigma)$.
The dependence of \eqref{eq:topping} on positive mean curvature alone is essential here.
Write $H_f$ for the scalar mean curvature of $f$ with respect to a local unit normal to $f$.
Away from the boundary, the positive mean-curvature contributions from the two local normal perturbations converge to $(H_f)_++(-H_f)_+=|\mathbf H_f|$, so the original curvature is not counted twice; at each point, at most one orientation contributes.
This avoids the factor-two loss of an estimate involving absolute mean curvature, as emphasized in \cite{Topping2026}.

The additional difficulty in the immersed case is that the doubled surface may still have self-intersections, as in \Cref{fig:doubling-sequence}, so one cannot directly apply \eqref{eq:topping} to the doubled surface.
However, after a generic perturbation, the doubled surface has only normal crossings (transverse self-intersections), and the exposed crossing parts are all ``valley'' folds.
Geometrically, these valley folds carry ``singular negative mean curvature'', and hence they are expected not to contribute to the positive mean-curvature energy.

Our key step makes this observation precise, without applying mean-convex hull regularity directly to a self-intersecting image.
Instead, we smooth the distance function to the image of an immersed solid bounded by the doubled surface, and take a suitable level set which is a smooth closed embedded surface enclosing the original immersed solid.
Applying Topping's embedded estimate to the resulting smooth outer approximations, we obtain an immersed-solid generalization of estimate \eqref{eq:topping}, see \Cref{thm:solid}.
Combining this estimate with the doubling trick proves \Cref{thm:main}.

In \Cref{sec:solids} we prove the immersed-solid estimate, and in \Cref{sec:tubes} we use the doubling trick to prove \Cref{thm:main}.

\begin{acknowledgements}
	The author would like to thank Peter Topping for helpful discussions.
	This work is supported by JSPS KAKENHI Grant Numbers JP23H00085 and JP24K00532.
\end{acknowledgements}

\section{An outward estimate for immersed solids}\label{sec:solids}

For a nonempty compact set $A\subset\R^3$, set
\[
\mathcal W(A)\coloneqq \frac12\int_{\Sph}W_A(e)\,d\HH^2(e).
\]
Thus $\mathcal W(A)$ is $2\pi$ times the mean width of $A$.

An \emph{immersed solid} is a smooth map $F:\Omega\to\R^3$ from a compact smooth three-manifold $\Omega$ with boundary $\partial\Omega$ such that $dF$ has rank three everywhere, including on $\partial \Omega$.
Equivalently, $F$ is a local diffeomorphism in the interior and extends locally across each boundary point as a diffeomorphism; the image of a boundary neighborhood in $\Omega$ locally occupies one side of its smooth boundary surface.
The restriction $f\coloneqq F|_{\partial \Omega}$ has a well-defined outward unit normal $\nu$ from these local domains.
We write $H_F$ for the scalar mean curvature of $f$ with respect to the outward unit normal $\nu$.
For a smooth boundary $M$ of a Euclidean domain, $H_M$ denotes its scalar mean curvature with the same outward convention.

The goal of this section is to prove the following extension of \eqref{eq:topping}.

\begin{theorem}\label[theorem]{thm:solid}
	Let $F:\Omega\to\R^3$ be a compact connected immersed solid, with possibly disconnected boundary $\partial\Omega$.
    Then 
	\begin{equation}\label{eq:solid}
		\mathcal W(F(\Omega))\leq\int_{\partial \Omega}(H_F)_+\,d\mu_f.
	\end{equation}
\end{theorem}

We first prove the key regularization property needed for \Cref{thm:solid}.
Recall that $f:\partial\Omega\to\R^3$ has \emph{normal crossings} if, for every $x\in f(\partial\Omega)$,
\[
	f^{-1}(x)=\{p_1,\ldots,p_k\},\qquad 1\leq k\leq3,
\]
and unit normals $\nu_i\perp df_{p_i}(T_{p_i}\partial\Omega)$ satisfy
\[
	\operatorname{rank}[\nu_1\ \cdots\ \nu_k]=k.
\]
The condition is independent of the signs of the $\nu_i$.
See \cite{GolubitskyGuillemin1973}*{Chapter~III, Definition~3.1 and Lemma~3.6} for normal crossings and the corresponding tangent-space condition.
Throughout this section, for $A\subset\R^3$, $x\in\R^3$, and $r>0$, we write
\[
	\begin{aligned}
		N_r(A)&\coloneqq \{x\in\R^3:\dist(x,A)<r\},\qquad \dist(x,A)\coloneqq \inf_{y\in A}|x-y|,\\
		B_r(x)&\coloneqq \{y\in\R^3:|y-x|<r\}.
	\end{aligned}
\]

\begin{lemma}\label[lemma]{lem:regularization}
	Let $F:\Omega\to\R^3$ be an immersed solid whose boundary immersion $f:\partial\Omega\to\R^3$ has only normal crossings.
	Then there are bounded open sets $U_\varepsilon\subset\R^3$ with smooth boundary $\partial U_\varepsilon$ such that $F(\Omega)\subset U_\varepsilon$ and
	\begin{equation}\label{eq:regularization}
		\limsup_{\varepsilon\downarrow0}
		\int_{\partial U_\varepsilon}(H_{\partial U_\varepsilon})_+\,d\HH^2
		\leq\int_{\partial \Omega}(H_F)_+\,d\mu_f,
	\end{equation}
	where the curvatures on $\partial U_\varepsilon$ are also defined through the outward normal.
\end{lemma}

\begin{proof}
    We will construct $U_\varepsilon$ by slightly thickening $F(\Omega)$ and then smoothing the singularities in $\partial F(\Omega)$.
    
    \emph{Step 1: The boundary self-intersection set.}
	We define the set of ordered pairs of distinct points with a common image by
	\[
		\mathcal D_2(f)\coloneqq\{(p,q):p,q\in\partial\Omega,\ p\ne q,\ f(p)=f(q)\}.
	\]
	Put $\Delta\coloneqq\{(p,p):p\in\partial\Omega\}$ and define $G:(\partial\Omega\times\partial\Omega)\setminus\Delta\to\R^3$ by $G(p,q)\coloneqq f(p)-f(q)$.
	The normal-crossing condition implies that the planes $df_p(T_p\partial\Omega)$ and $df_q(T_q\partial\Omega)$ span $\R^3$ for every $(p,q)\in\mathcal D_2(f)$, so $0$ is a regular value of $G$.
	Thus $\mathcal D_2(f)=G^{-1}(0)$ is a smooth one-dimensional manifold (without boundary), possibly empty.

	Let $d_{\partial\Omega}$ be the Riemannian distance induced by $f$, with distance $+\infty$ between distinct connected components.
	Compactness of $\partial\Omega$ together with local injectivity of $f$ gives $a>0$ such that $d_{\partial\Omega}(p,q)\geq a$ holds for all $(p,q)\in\mathcal D_2(f)$.
	Thus $\mathcal D_2(f)$ is closed in a compact subset of $(\partial\Omega\times\partial\Omega)\setminus\Delta$, hence compact.
    
	Let $\pi_1(p,q)\coloneqq p$ and define the common-image set by
	\[
		D\coloneqq\{f(p):(p,q)\in\mathcal D_2(f)\}
		=(f\circ\pi_1)(\mathcal D_2(f))\subset\R^3.
	\]
	The restriction of $f\circ\pi_1$ to $\mathcal D_2(f)$ is an immersion.
	Indeed, for $(p,q)\in\mathcal D_2(f)$ and $(v,w)\in T_{(p,q)}\mathcal D_2(f)$, we have
	\[
		df_p(v)=df_q(w),\qquad d(f\circ\pi_1)_{(p,q)}(v,w)=df_p(v).
	\]
	If the latter derivative vanishes, injectivity of $df_p$ gives $v=0$, and then injectivity of $df_q$ gives $w=0$.
	Thus the differential has zero kernel.
    
	Consequently, $D$ is a finite union of smooth immersed closed curves.
	(The same argument also shows that the set of triple points is finite, although this fact will not be used below.)
    
	When $D\ne\varnothing$, as the total length of $D$ is finite, there are finitely many points $x_1,\ldots,x_{J_r}\in D$ such that
	\begin{equation}\label{eq:double-cover}
		D\subset\bigcup_{j=1}^{J_r}B_r(x_j),
		\qquad J_r\leq\frac{C_D}{r}\quad(0<r\leq1),
	\end{equation}
    where $C_D>0$ is independent of $r$.
	If $D=\varnothing$, set $J_r\coloneqq0$ and $C_D\coloneqq1$.

	\emph{Step 2: Local boundary pieces and signed distances.}
	Let $q\in\partial(F(\Omega))$ be arbitrary.
	Since $F$ is a local diffeomorphism on $\operatorname{int}\Omega$, the set $F(\operatorname{int}\Omega)$ is open and hence contained in $\operatorname{int}(F(\Omega))$.
	Thus $F^{-1}(q)\subset\partial\Omega$.
	By the normal-crossing assumption,
	\[
		F^{-1}(q)=f^{-1}(q)=\{p_1,\ldots,p_k\}\subset\partial\Omega,
		\qquad 1\leq k\leq3.
	\]
	Choose pairwise disjoint neighborhoods $V_i$ of $p_i$ in $\Omega$ such that $F|_{V_i}:V_i\to F(V_i)$ is a diffeomorphism onto a region with smooth boundary near $q$.
	By compactness, we can choose a sufficiently small ball $B_q\subset\R^3$ centered at $q$ such that
	\[
		F^{-1}(B_q)\subset\bigcup_{i=1}^kV_i.
	\]
	There are thus smooth signed-distance functions $\phi_1,\ldots,\phi_k$ on $B_q$ to the local sheets $f(V_i\cap\partial\Omega)$, with $\phi_i>0$ on $B_q\setminus F(V_i)$, such that
	\begin{equation}\label{eq:branches}
		F(\Omega)\cap B_q=\bigcup_{i=1}^k\{x\in B_q:\phi_i(x)\leq0\}.
	\end{equation}
	Further shrinking $B_q$ if necessary, we may assume that nearest-point projections from $B_q$ onto each $f(V_i\cap\partial\Omega)$ are well defined and nearest points in $F(\Omega)$ to points of $B_q$ lie in $\bigcup_i F(V_i)$.
    Write $d_{F(\Omega)}(x)\coloneqq\operatorname{dist}(x,F(\Omega))$ for $x\in\R^3$.
	Then
	\begin{equation}\label{eq:distance-branches}
		d_{F(\Omega)}(x)=\min_i\phi_i(x)\quad(x\in B_q\setminus F(\Omega)).
	\end{equation}
	The linear independence of the normals gives
	\[
		0\notin\operatorname{conv}\{\nabla\phi_1(q),\ldots,\nabla\phi_k(q)\}.
	\]
	By the strict separation theorem, there exist $v\in\Sph$ and $c>0$ such that $\langle v,\nabla\phi_i(q)\rangle\geq2c$ for every $i$.
	By continuity, after shrinking $B_q$, we have
	\begin{equation}\label{eq:local-directional-bound}
		\langle v,\nabla\phi_i(x)\rangle\geq c\qquad(x\in B_q,\ 1\leq i\leq k).
	\end{equation}
    \Cref{fig:local-regularization}(a) illustrates the local union in \eqref{eq:branches} for $k=2$.

	As $q$ varies over $\partial(F(\Omega))$, we choose $r_q>0$ with $\overline{B_{4r_q}(q)}\subset B_q$.
	Compactness gives a finite subcover $B_{r_j}(a_j)$, $1\leq j\leq J$, of $\partial(F(\Omega))$, where $r_j\coloneqq r_{a_j}$.
	For each $j$, let $V_{j,1},\ldots,V_{j,k_j}$ denote the neighborhoods constructed at $q=a_j$ and $v_j$ the corresponding unit vector.
	Let $\phi_{j,i}$ be the signed distance to $f(V_{j,i}\cap\partial\Omega)$, with $\phi_{j,i}>0$ on $B_{4r_j}(a_j)\setminus F(V_{j,i})$.
	These functions are smooth near $\overline{B_{4r_j}(a_j)}$ and satisfy \eqref{eq:branches}, \eqref{eq:distance-branches}, and \eqref{eq:local-directional-bound} there.
	Since the covering is finite, there are $C>0$ and $c>0$, independent of $j$ and $i$, such that
	\begin{equation}\label{eq:cone}
		\begin{gathered}
			\|D^2\phi_{j,i}(x)\|\leq C,
			\qquad \langle v_j,\nabla\phi_{j,i}(x)\rangle\geq c>0,\\
			x\in B_{4r_j}(a_j),\qquad 1\leq j\leq J,\quad1\leq i\leq k_j.
		\end{gathered}
	\end{equation}

	\emph{Step 3: Exterior level sets of the smoothed distance.}
    Given the above setup, we now smooth the boundary $\partial(F(\Omega))$ by mollifying the distance function $d_{F(\Omega)}$.
	Choose
	\[
		\rho\in C_c^\infty(B_1(0)),\qquad \rho\geq0,\qquad
		\int_{\R^3}\rho\,dx=1,\qquad \rho_\delta(x)\coloneqq\delta^{-3}\rho(x/\delta),
	\]
	and for $0<\varepsilon<1$, set
	\begin{equation}\label{eq:smoothing}
		\delta\coloneqq\varepsilon^2,
		\qquad g_\varepsilon\coloneqq d_{F(\Omega)}*\rho_\delta,
		\qquad U_\varepsilon\coloneqq\{g_\varepsilon<\varepsilon\},
		\qquad M_\varepsilon\coloneqq\{g_\varepsilon=\varepsilon\}.
	\end{equation}
	\Cref{fig:local-regularization}(b) illustrates the smoothed exterior level $M_\varepsilon$.

	\begin{figure}[htbp]
		\centering
		\begin{minipage}[t]{.48\linewidth}
			\centering
			\includegraphics[width=\linewidth]{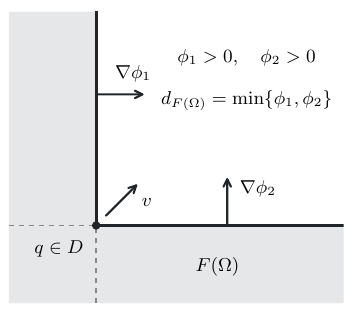}
			\par\smallskip
			\normalfont (a) The local union of boundary pieces
		\end{minipage}\hfill
		\begin{minipage}[t]{.48\linewidth}
			\centering
			\includegraphics[width=\linewidth]{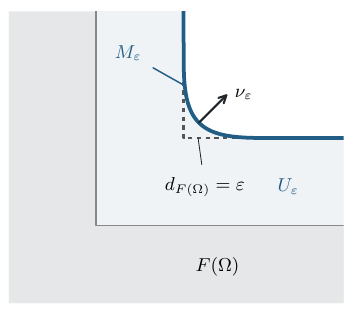}
			\par\smallskip
			\normalfont (b) Smoothed exterior level
		\end{minipage}
		\caption{Local regularization near a crossing point $q\in D$ where two boundary pieces meet transversely, shown in a cross-section perpendicular to the intersection curve.
		In (a), the gray region is $F(\Omega)$, while $d_{F(\Omega)}=\min\{\phi_1,\phi_2\}>0$ holds in the exterior white region.
		In (b), gray represents the original $F(\Omega)$, and its union with blue represents the smoothed region $U_\varepsilon$.
		The dashed curve is $\{d_{F(\Omega)}=\varepsilon\}$, and the blue curve is the smoothed level $M_\varepsilon$.}
		\label{fig:local-regularization}
	\end{figure}

	Since $d_{F(\Omega)}$ is $1$-Lipschitz, we have $|g_\varepsilon-d_{F(\Omega)}|\leq\delta$, $|\nabla g_\varepsilon|\leq1$, and
	\[
		F(\Omega)\subset U_\varepsilon\subset N_{\varepsilon+\delta}(F(\Omega)).
	\]
	For $x\in M_\varepsilon$, we have $0<\varepsilon-\delta\leq d_{F(\Omega)}(x)\leq\varepsilon+\delta$; the positivity of $d_{F(\Omega)}$ yields $\operatorname{dist}(x,\partial(F(\Omega)))=d_{F(\Omega)}(x)$.
	We choose $0<\varepsilon_0<1/4$ with $\varepsilon_0+\varepsilon_0^2<\min_jr_j$.
	Since $\{B_{r_j}(a_j)\}$ covers $\partial(F(\Omega))$, the bound $\operatorname{dist}(x,\partial(F(\Omega)))\leq\varepsilon+\varepsilon^2$ implies that every $x\in M_\varepsilon$ belongs to some $B_{2r_j}(a_j)$ for $0<\varepsilon<\varepsilon_0$.
    
	For the rest of the proof, we assume $0<\varepsilon<\varepsilon_0$.
    
	For each $j$ and $x\in M_\varepsilon\cap B_{2r_j}(a_j)$, the inequalities $2\delta=2\varepsilon^2<\varepsilon_0+\varepsilon_0^2<r_j$ give $B_{2\delta}(x)\subset B_{3r_j}(a_j)$.
	We have $d_{F(\Omega)}\geq\varepsilon-3\delta>0$ on $B_{2\delta}(x)$, so \eqref{eq:distance-branches} applies and $\phi_{j,i}>0$ for $1\leq i\leq k_j$.
	At almost every $y\in B_{2\delta}(x)$, the function $d_{F(\Omega)}$ is differentiable; choosing $i$ with $d_{F(\Omega)}(y)=\phi_{j,i}(y)$ then gives $\nabla d_{F(\Omega)}(y)=\nabla\phi_{j,i}(y)$.
	Hence convolution and \eqref{eq:cone} give
	\begin{equation}\label{eq:directional-bound}
		\partial_{v_j}g_\varepsilon(x)
		=\int_{B_\delta(0)}\partial_{v_j}d_{F(\Omega)}(x-y)\rho_\delta(y)\,dy\geq c.
	\end{equation}
	Since $v_j$ is a unit vector and $x\in M_\varepsilon$ was arbitrary, we obtain
	\begin{equation}\label{eq:gradient-lower-bound}
		|\nabla g_\varepsilon(x)|\geq c\qquad(x\in M_\varepsilon).
	\end{equation}
    In particular, the implicit function theorem shows that $M_\varepsilon$ is smooth and
	\[
		M_\varepsilon=\partial U_\varepsilon.
	\]
    
	Moreover, each function $y\mapsto\phi_{j,i}(y)-C|y|^2/2$ is concave on $B_{4r_j}(a_j)$, and hence their minimum is concave.
	On $B_{2\delta}(x)$, where $d_{F(\Omega)}=\min_i\phi_{j,i}$, this minimum is
	\[
		\min_{1\leq i\leq k_j}\bigl(\phi_{j,i}(y)-C|y|^2/2\bigr)=d_{F(\Omega)}(y)-C|y|^2/2.
	\]
	Thus $D^2d_{F(\Omega)}\leq CI$ distributionally on $B_{2\delta}(x)$, and, after convolution, $D^2g_\varepsilon(x)\leq CI$ pointwise, where $I$ is the identity matrix.
	Set $\nu_\varepsilon\coloneqq\nabla g_\varepsilon/|\nabla g_\varepsilon|$.
	At each $x\in M_\varepsilon$, the matrix $I-\nu_\varepsilon\otimes\nu_\varepsilon$ is the orthogonal projection onto $T_xM_\varepsilon$.
	Choosing an orthonormal basis $e_1,e_2$ of the tangent plane $T_xM_\varepsilon$, the pointwise Hessian bound gives
	\[
		\tr((I-\nu_\varepsilon\otimes\nu_\varepsilon)D^2g_\varepsilon)
		=\sum_{i=1}^2\langle D^2g_\varepsilon e_i,e_i\rangle
		\leq\sum_{i=1}^2 C|e_i|^2=2C.
	\]
	The level-set formula for mean curvature and \eqref{eq:gradient-lower-bound} now yield
	\begin{equation}\label{eq:level-bounds}
		H_{M_\varepsilon}
		=\frac{\tr((I-\nu_\varepsilon\otimes\nu_\varepsilon)D^2g_\varepsilon)}{2|\nabla g_\varepsilon|}
		\leq\frac Cc.
	\end{equation}

	\emph{Step 4: Area estimates near crossings.}
    In the remaining steps, we estimate the energy of the smoothed surfaces $M_\varepsilon$.
    We first show that the energy near the crossing part $D$ is negligible.
    
	Set $E_j\coloneqq M_\varepsilon\cap B_{2r_j}(a_j)$ and let
	\[
		\pi_{v_j}:\R^3\to v_j^\perp,\qquad
		\pi_{v_j}(x)\coloneqq x-\langle x,v_j\rangle v_j.
	\]
	Then $\pi_{v_j}|_{E_j}$ is injective.
	Indeed, we fix $z\in v_j^\perp$ and consider $h(t)\coloneqq g_\varepsilon(z+tv_j)-\varepsilon$ on the interval $I_z\coloneqq\{t\in\R:z+tv_j\in B_{2r_j}(a_j)\}$.
	By \eqref{eq:directional-bound}, $h'(t)\geq c>0$ at every zero of $h$, so each zero is a strict crossing from negative to positive values.
	Thus $h$ has at most one zero, and hence no distinct pair of points in $E_j$ is projected to the same $z$.

	On each $E_j$, the area-change factor of the projection $\pi_{v_j}|_{M_\varepsilon}$ is $|\langle\nu_\varepsilon,v_j\rangle|$, which is at least $c$ by \eqref{eq:directional-bound} and $|\nabla g_\varepsilon|\leq1$.
	The area formula consequently gives
	\[
		\HH^2(E_j\cap B_r(x))
		\leq c^{-1}\HH^2(\pi_{v_j}(E_j\cap B_r(x)))
		\leq c^{-1}\pi r^2.
	\]
	Summing over $1\leq j\leq J$ yields, for $x\in\R^3$, $r>0$, and $0<\varepsilon<\varepsilon_0$,
	\begin{equation}\label{eq:local-area}
		\HH^2(M_\varepsilon\cap B_r(x))\leq C_1r^2,\qquad C_1\coloneqq J\pi/c.
	\end{equation}
	The cover in \eqref{eq:double-cover} gives $N_r(D)\subset\bigcup_{j=1}^{J_r}B_{2r}(x_j)$ with $J_r\leq C_D/r$.
	Using this inclusion, \eqref{eq:local-area}, and \eqref{eq:level-bounds}, set $C_2\coloneqq 4C_1C_D$ and $C_3\coloneqq (C/c)C_2$, which are independent of $\varepsilon$ and $r$.
	Then, for $0<\varepsilon<\varepsilon_0$ and $0<r\leq1$,
	\begin{equation}\label{eq:crossing-area}
		\HH^2(M_\varepsilon\cap N_r(D))\leq C_2r,
		\qquad
		\int_{M_\varepsilon\cap N_r(D)}(H_{M_\varepsilon})_+\,d\HH^2\leq C_3r.
	\end{equation}

	\emph{Step 5: Boundary points with one preimage.}
	We then estimate the energy contribution of the non-crossing part.
    
    At $q\in\partial(F(\Omega))\setminus D$, the set $f^{-1}(q)$ consists of one point, so $k=1$ in \eqref{eq:branches} and \eqref{eq:distance-branches}.
	Since $\varepsilon-2\delta>0$, for $x\in M_\varepsilon$ sufficiently near $q$, the convolution defining $g_\varepsilon(x)$ in \eqref{eq:smoothing} only samples points outside $F(\Omega)$, where $d_{F(\Omega)}=\phi_1$.
	Thus $g_\varepsilon=\phi_1*\rho_\delta$ in a neighborhood of $M_\varepsilon$ near $q$.
	Since $\phi_1$ is smooth, $\phi_1*\rho_\delta\to\phi_1$ in $C^2$ locally as $\varepsilon\downarrow0$.
	As $|\nabla\phi_1|=1$, the implicit function theorem gives $C^2$ convergence of the nearby part of $M_\varepsilon$ to $\partial(F(\Omega))$ near $q$, with the same outward orientation.
    
	For a nonnegative $\psi\in C_c^\infty(\R^3\setminus D)$, take a finite partition of unity over such neighborhoods covering the compact set $\partial(F(\Omega))\cap\operatorname{supp}\psi$.
	Since $f$ is injective on $f^{-1}(\partial(F(\Omega))\setminus D)$, the $C^2$ convergence and change of variables under $f$ give
	\begin{equation}\label{eq:exposed-limit}
		\begin{aligned}
		\lim_{\varepsilon\downarrow0}
		\int_{M_\varepsilon}(H_{M_\varepsilon})_+\psi\,d\HH^2
		&=\int_{f^{-1}(\partial(F(\Omega))\setminus D)}
		(H_F)_+(\psi\circ F)\,d\mu_f\\
		&\leq\int_{\partial\Omega}(H_F)_+(\psi\circ F)\,d\mu_f.
		\end{aligned}
	\end{equation}
	The inequality follows from $f^{-1}(\partial(F(\Omega))\setminus D)\subset\partial\Omega$ and $(H_F)_+(\psi\circ F)\geq0$.
	
	\emph{Step 6: Conclusion.}
	Choose $R>0$ so that $M_\varepsilon\subset B_R(0)$ for all small $\varepsilon$, and choose $\chi_r\in C_c^\infty(B_{2R}(0))$ such that
	\[
		0\leq\chi_r\leq1,\qquad
		\chi_r=0\text{ on }N_r(D),\qquad
		\chi_r=1\text{ on }B_R(0)\setminus N_{2r}(D).
	\]
	By \eqref{eq:crossing-area}, for $0<r\leq1/2$,
	\[
		\int_{M_\varepsilon}(H_{M_\varepsilon})_+(1-\chi_r)\,d\HH^2
		\leq\int_{M_\varepsilon\cap N_{2r}(D)}(H_{M_\varepsilon})_+\,d\HH^2
		\leq2C_3r.
	\]
	Applying \eqref{eq:exposed-limit} with $\psi=\chi_r$ therefore gives
	\[
		\limsup_{\varepsilon\downarrow0}
		\int_{M_\varepsilon}(H_{M_\varepsilon})_+\,d\HH^2
		\leq\int_{\partial \Omega}(H_F)_+\,d\mu_f+2C_3r.
	\]
	Letting $r\downarrow0$ proves \eqref{eq:regularization}.
	If $D$ is empty, the crossing estimates are unnecessary and \eqref{eq:exposed-limit} with a cutoff equal to one near $F(\Omega)$ gives the result directly.
\end{proof}

\begin{proof}[Proof of \Cref{thm:solid}]
	We first address the generic case where the boundary only has normal crossings, and then complete the general case by perturbations.
    
    \emph{Step 1: The normal-crossing case.}
	First suppose that $f$ has normal crossings and use $U_\varepsilon$ from \Cref{lem:regularization}.
	The connected set $F(\Omega)$ lies in one connected component $U_\varepsilon^0$ of $U_\varepsilon$.
	Let $O_\varepsilon$ be the unbounded connected component of $\R^3\setminus\overline{U_\varepsilon^0}$.
	Choose a connected component $M_\varepsilon^0$ of $\partial U_\varepsilon^0$ such that $M_\varepsilon^0\cap\partial O_\varepsilon\ne\varnothing$.
	By the Jordan--Brouwer separation theorem and the connectedness of $U_\varepsilon^0$, the bounded open region $V_\varepsilon^0$ enclosed by $M_\varepsilon^0$ contains $U_\varepsilon^0$, so
	\[
		\partial V_\varepsilon^0=M_\varepsilon^0,
		\qquad F(\Omega)\subset U_\varepsilon^0\subset V_\varepsilon^0.
	\]

	The outward normal of $M_\varepsilon^0$ agrees with that of $\partial U_\varepsilon$ on $M_\varepsilon^0 \subset \partial U_\varepsilon$.
	Observing that $W_{F(\Omega)}\leq W_{\overline{V_\varepsilon^0}}= W_{M_\varepsilon^0}$, and applying \eqref{eq:topping} to $M_\varepsilon^0$ and discarding the other nonnegative boundary contributions, we obtain
	\[
		\mathcal W(F(\Omega)) \leq \mathcal W(M_\varepsilon^0)
		\leq\int_{M_\varepsilon^0}(H_{M_\varepsilon^0})_+\,d\HH^2
		\leq\int_{\partial U_\varepsilon}(H_{\partial U_\varepsilon})_+\,d\HH^2.
	\]
	The limit in \eqref{eq:regularization} proves \eqref{eq:solid} in this case.

	\emph{Step 2: Approximation of a general immersed solid.}
	For general $F$, a standard perturbation argument (see, e.g., \cite{GolubitskyGuillemin1973}*{Chapter~III, Corollary~3.3}) gives smooth normal-crossing immersions $f_j:\partial\Omega\to\R^3$ converging to $f$ in $C^2$.
	Choose a smooth embedding $c:\partial\Omega\times[0,a)\to\Omega$ with $c(p,0)=p$, whose image is a neighborhood of $\partial\Omega$, and extend $f_j$ to $F_j$ by
	\[
		F_j(c(p,t))\coloneqq F(c(p,t))+\beta(t)\bigl(f_j(p)-f(p)\bigr),
	\]
	where $\beta\in C_c^\infty([0,a))$ and $\beta(0)=1$, and put $F_j\coloneqq F$ outside $c(\partial\Omega\times\operatorname{supp}\beta)$.
	In particular, $F_j|_{\partial\Omega}=f_j$.
	Then $\|F_j-F\|_{C^2(\Omega)}\rightarrow0$ as $j\to\infty$, while $\operatorname{rank}dF_j(p)=3$ for all $p\in\Omega$ and all sufficiently large $j$ by compactness.
	By the $C^2$ convergence, $\mathcal W(F_j(\Omega))\to\mathcal W(F(\Omega))$ and $\int_{\partial\Omega}(H_{F_j})_+\,d\mu_{f_j}\to\int_{\partial\Omega}(H_F)_+\,d\mu_f$.
	Passing to the limit in \eqref{eq:solid} for $F_j$ proves the theorem.
\end{proof}

\section{The doubling trick}\label{sec:tubes}

Now we complete the proof of \Cref{thm:main} via the doubling construction as follows.

\begin{proposition}\label[proposition]{prop:tube}
	Let $f:\Sigma\to\R^3$ be an immersion of a compact connected surface $\Sigma$, possibly with boundary $\partial\Sigma$.
	Then there are compact connected immersed solids $F_\varepsilon:\Omega_\varepsilon\to\R^3$ such that $f(\Sigma)\subset F_\varepsilon(\Omega_\varepsilon)$ and
	\begin{equation}\label{eq:tube-energy}
		\lim_{\varepsilon\downarrow0}
		\int_{\partial \Omega_\varepsilon}(H_{F_\varepsilon})_+\,d\mu_{f_\varepsilon}
		=\int_\Sigma|\mathbf H_f|\,d\mu_f+\frac\pi2\operatorname{Length}(\partial\Sigma),
	\end{equation}
    where $f_\varepsilon\coloneqq F_\varepsilon|_{\partial\Omega_\varepsilon}:\partial\Omega_\varepsilon\to\R^3$.
\end{proposition}

\begin{proof}
	The general idea follows the doubling construction of \cite{Miura2022}*{Section~2, Proof of Proposition~2.1}, with the positive-part curvature energy as in \cite{Topping2026}*{Proof of Theorem~1.8}.
    Here we describe the construction in terms of the normal bundle of an immersion extending $f$ across $\partial\Sigma$ so that self-intersections cause no ambiguity; see \Cref{fig:doubling-sequence}.

	\emph{Step 1: Construction.}
	Extend $f$ across a small boundary collar to a compact immersed surface $\widetilde f:\widetilde\Sigma\to\R^3$ with $\Sigma\subset\operatorname{int}\widetilde\Sigma$; if $\partial\Sigma=\varnothing$, take $\widetilde\Sigma\coloneqq\Sigma$ and $\widetilde f\coloneqq f$.
	Write the normal line bundle of $\widetilde f$ as
	\[
		N\widetilde\Sigma\coloneqq\{(p,v)\in\widetilde\Sigma\times\R^3:v\perp d\widetilde f(T_p\widetilde\Sigma)\},
	\]
	and the restriction of this bundle to $\Sigma$ as $N\Sigma$.
	Define
	\[
		F:N\widetilde\Sigma\to\R^3,\qquad F(p,v)\coloneqq\widetilde f(p)+v.
	\]
	Since $dF$ is invertible along the zero section, compactness gives $\tau_0>0$ such that $F$ is an immersion on $\{(p,v)\in N\widetilde\Sigma:|v|<\tau_0\}$.

	Choose pairwise disjoint tubular neighborhoods $U_\Gamma\subset\widetilde\Sigma$ of the connected components $\Gamma$ of $\partial\Sigma$.
	Set $L_\Gamma\coloneqq\operatorname{Length}(\Gamma)$ and, for sufficiently small $r_0>0$, use Fermi coordinates $(u,r)\in(\R/L_\Gamma\mathbf Z)\times(-r_0,r_0)$ on each $U_\Gamma$, after shrinking $U_\Gamma$ if necessary, for the metric induced by $\widetilde f$, where $u$ is arclength along $\Gamma$, and $r$ is signed geodesic distance from $\Gamma$, positive inward, so that $\Sigma$ corresponds to $r\geq0$.
	Each $U_\Gamma$ is an annulus and hence admits a smooth unit normal; fix one choice, denoted by $\nu(u,r)$, on each $U_\Gamma$.
	Fix a smooth unit-speed convex arc $\gamma:[0,\ell_0]\to\R^2$, $\gamma(s)\coloneqq(a(s),b(s))$, from $(0,1)$ to $(0,-1)$, lying in $a\leq0$, and joining smoothly to the rays $b=\pm1$, $a\geq0$.
	Its signed curvature $\kappa\coloneqq a'b''-b'a''$ satisfies
	\begin{equation}\label{eq:profile-turning}
		\kappa\geq0,\qquad \int_0^{\ell_0}\kappa(s)\,ds=\pi.
	\end{equation}
	Let $P$ be the closed convex region bounded by this arc $\gamma$ and the two rays $b=\pm1$, $a\geq0$, and choose $C>0$ with $P\subset[-C,\infty)\times[-1,1]$.
	Choose $0<\varepsilon<\tau_0$, with $C\varepsilon<r_0$ when $\partial\Sigma\ne\varnothing$.
	For each boundary component $\Gamma$, set
	\[
        E_{\varepsilon,\Gamma}\coloneqq\bigl\{((u,r),t\nu(u,r)):\ u\in\R/L_\Gamma\mathbf Z,\ -r_0<r<0,\ (r/\varepsilon,t/\varepsilon)\in P\bigr\}\subset N\widetilde\Sigma.
	\]
	Then we define
	\[
		\Omega_\varepsilon\coloneqq\{(p,v)\in N\Sigma:|v|\leq\varepsilon\}
		\cup\bigcup_{\Gamma\subset\partial\Sigma}E_{\varepsilon,\Gamma}.
	\]
	See \Cref{fig:collar-section} for a cross-section.
	\begin{figure}[htbp]
		\centering
		\includegraphics[width=.48\linewidth]{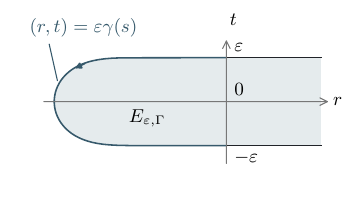}
		\caption{Cross-section of $\Omega_\varepsilon$ near a point of $\Gamma\subset\partial\Sigma$, in the $(r,t)$ coordinates at fixed $u$.}
		\label{fig:collar-section}
	\end{figure}

	The smooth joining of the profile and the bound $r\geq-C\varepsilon>-r_0$ show that $\Omega_\varepsilon$ is a compact smooth three-manifold with boundary.
	For $(p,v)\in\Omega_\varepsilon\cap N\Sigma$, the path $\lambda\mapsto(p,(1-\lambda)v)$, $0\leq\lambda\leq1$, lies in $\Omega_\varepsilon$ and joins $(p,v)$ to $(p,0)$.
	For $((u,r),t\nu(u,r))\in E_{\varepsilon,\Gamma}$, the path
	\[
		\lambda\mapsto\bigl((u,(1-\lambda)r),(1-\lambda)t\nu(u,(1-\lambda)r)\bigr),
		\qquad 0\leq\lambda\leq1,
	\]
	lies in $\Omega_\varepsilon$ and ends at $((u,0),0)$, since $P$ is convex and contains $(0,0)$.
	Thus every point is joined to the zero section over the connected surface $\Sigma$, so $\Omega_\varepsilon$ is connected.
	The restriction $F_\varepsilon\coloneqq F|_{\Omega_\varepsilon}:\Omega_\varepsilon\to\R^3$ is an immersed solid and satisfies $f(\Sigma)\subset F_\varepsilon(\Omega_\varepsilon)$.
	When $\partial\Sigma=\varnothing$, the collar union is empty and the same conclusions hold.

	\emph{Step 2: Energy convergence.}
	Define
	\[
		S_\varepsilon\coloneqq\{(p,v)\in N\Sigma:|v|=\varepsilon\},
	\]
	and, for each boundary component $\Gamma$,
	\[
		C_{\varepsilon,\Gamma}\coloneqq\bigl\{((u,\varepsilon a(s)),\varepsilon b(s)\nu(u,\varepsilon a(s))):u\in\R/L_\Gamma\mathbf Z,\quad 0\leq s\leq\ell_0\bigr\}.
	\]
	Then the boundary of $\Omega_\varepsilon$ is given by
	\begin{equation}\label{eq:tube-boundary-decomposition}
		\partial\Omega_\varepsilon=S_\varepsilon\cup\bigcup_{\Gamma\subset\partial\Sigma}C_{\varepsilon,\Gamma}.
	\end{equation}
	The caps $C_{\varepsilon,\Gamma}$ are pairwise disjoint, and each intersection $S_\varepsilon\cap C_{\varepsilon,\Gamma}$ has zero area.
	Around each point of $\Sigma$, choose a smooth unit normal field $\nu$ to $f$.
	Locally, the restriction of $f_\varepsilon$ to $S_\varepsilon$ is represented by two sheets $f\pm\varepsilon\nu$, with opposite outward unit normals $\pm\nu$.
	Writing $H_f$ for the scalar mean curvature of $f$ with respect to $\nu$, as $\varepsilon\downarrow0$ the scalar mean curvatures converge to $H_f$ and $-H_f$, respectively, and the area densities relative to $d\mu_f$ converge to $1$, uniformly on compact coordinate patches.
	A finite cover of $\Sigma$ makes these convergences uniform for all $p\in\Sigma$.
	Thus, since $(H_f)_++(-H_f)_+=|\mathbf H_f|$,
	\begin{equation}\label{eq:normal-face-limit}
		\lim_{\varepsilon\downarrow0}\int_{S_\varepsilon}(H_{F_\varepsilon})_+\,d\mu_{f_\varepsilon}
		=\int_\Sigma|\mathbf H_f|\,d\mu_f.
	\end{equation}
	Note that $\nu$ is chosen only locally; no global unit normal field on $\Sigma$ is required.

	For each boundary component $\Gamma$, we define $X_\varepsilon:(\R/L_\Gamma\mathbf Z)\times[0,\ell_0]\to\R^3$ by
	\[
		X_\varepsilon(u,s)\coloneqq\widetilde f(u,\varepsilon a(s))+\varepsilon b(s)\nu(u,\varepsilon a(s)).
	\]
	This parametrizes the image of $C_{\varepsilon,\Gamma}$ under $F_\varepsilon$.
	Writing $g_{ij}\coloneqq\langle\partial_iX_\varepsilon,\partial_jX_\varepsilon\rangle$ for $i,j\in\{u,s\}$, the same local calculation as in \cite{Miura2022}*{Section~2, Proof of Proposition~2.1}, performed here in Fermi coordinates, gives $g_{uu}=1+O(\varepsilon)$, $g_{us}=O(\varepsilon^2)$, and $g_{ss}=\varepsilon^2(1+O(\varepsilon))$.
	The outward unit normal satisfies $n_\varepsilon=b'\partial_r\widetilde f-a'\nu+O(\varepsilon)$, where $\partial_r\widetilde f$ and $\nu$ are evaluated at $(u,\varepsilon a(s))$.
	Direct differentiation gives $\langle\partial_{ss}X_\varepsilon,n_\varepsilon\rangle=-\varepsilon\kappa(s)+O(\varepsilon^2)$, while $\partial_{uu}X_\varepsilon=O(1)$ and $\partial_{us}X_\varepsilon=O(\varepsilon)$.
	The formulas for mean curvature and area therefore yield, as $\varepsilon\downarrow0$,
	\begin{equation}\label{eq:cap-limit}
		H_{F_\varepsilon}=\frac{\kappa(s)}{2\varepsilon}+O(1),
		\qquad d\mu_{f_\varepsilon}=\varepsilon(1+O(\varepsilon))\,du\,ds.
	\end{equation}
	The error bounds in \eqref{eq:cap-limit} are uniform for $(u,s)\in(\R/L_\Gamma\mathbf Z)\times[0,\ell_0]$ and all boundary components $\Gamma$, with constants depending only on $\widetilde f$, the fixed neighborhoods $U_\Gamma$, and $\gamma$.
	Since $\kappa\geq0$ and $t\mapsto t_+$ is $1$-Lipschitz, we also have $(H_{F_\varepsilon})_+=\kappa(s)/(2\varepsilon)+O(1)$ uniformly.
    Together with \eqref{eq:profile-turning} and
    \eqref{eq:cap-limit}, this yields
	\begin{equation}\label{eq:cap-energy}
		\int_{C_{\varepsilon,\Gamma}}(H_{F_\varepsilon})_+\,d\mu_{f_\varepsilon}
		=\frac12\int_0^{L_\Gamma}\!du\int_0^{\ell_0}\kappa(s)\,ds+O(\varepsilon)
		=\frac\pi2L_\Gamma+O(\varepsilon).
	\end{equation}
	Combining \eqref{eq:normal-face-limit} and \eqref{eq:cap-energy} with \eqref{eq:tube-boundary-decomposition} proves \eqref{eq:tube-energy}.
\end{proof}

\begin{proof}[Proof of \Cref{thm:main}]
	Let $F_\varepsilon:\Omega_\varepsilon\to\R^3$ be the immersed solids given by \Cref{prop:tube}.
	For each fixed small $\varepsilon$, \Cref{thm:solid} and $f(\Sigma)\subset F_\varepsilon(\Omega_\varepsilon)$ give
	\[
		\mathcal W(f(\Sigma))
		\leq\mathcal W(F_\varepsilon(\Omega_\varepsilon))
		\leq\int_{\partial \Omega_\varepsilon}(H_{F_\varepsilon})_+\,d\mu_{f_\varepsilon}.
	\]
	We let $\varepsilon\downarrow0$ in \eqref{eq:tube-energy} to obtain the first inequality in \eqref{eq:main}.

	Since $f$ is an immersion, $f(\Sigma)$ is not contained in a line.
	Taking the convex hull preserves both directional widths and diameter.
	The strict inequality in \eqref{eq:main} follows from the classical strict monotonicity of mean width under inclusion, applied to the convex hull of $f(\Sigma)$ and a segment realizing its diameter.
	See also \cite{Topping2026}*{Remark~1.6} for a direct proof.
\end{proof}

\begin{AI}
    ChatGPT was used interactively in the preparation and revision of this manuscript. The author independently reviewed and revised all mathematical statements, proofs, and exposition and takes full responsibility for the final manuscript.
\end{AI}

\bibliography{references}

@misc{Topping2026,
	author = {Topping, Peter M.},
	title = {Sharp diameter estimates via integral geometry},
	year = {2026},
	note = {Preprint, \href{https://arxiv.org/abs/2610.34367v1}{arXiv:2610.34367v1}},
	url = {https://arxiv.org/abs/2610.34367v1}
}

@article{Topping1998,
	author = {Topping, Peter M.},
	title = {Mean curvature flow and geometric inequalities},
	journal = {J. Reine Angew. Math.},
	volume = {503},
	year = {1998},
	pages = {47--61},
	doi = {10.1515/crll.1998.099},
	note = {DOI \href{https://doi.org/10.1515/crll.1998.099}{\nolinkurl{10.1515/crll.1998.099}}}
}

@article{Miura2022,
	author = {Miura, Tatsuya},
	title = {A diameter bound for compact surfaces and the {Plateau--Douglas} problem},
	journal = {Ann. Sc. Norm. Super. Pisa Cl. Sci. (5)},
	volume = {23},
	number = {4},
	year = {2022},
	pages = {1707--1721},
	doi = {10.2422/2036-2145.202011_006},
	note = {DOI \href{https://doi.org/10.2422/2036-2145.202011_006}{\nolinkurl{10.2422/2036-2145.202011_006}}}
}

@article{Spadaro2020,
	author = {Spadaro, Emanuele},
	title = {Mean-convex sets and minimal barriers},
	journal = {Matematiche (Catania)},
	volume = {75},
	number = {1},
	year = {2020},
	pages = {353--375},
	doi = {10.4418/2020.75.1.16},
	note = {DOI \href{https://doi.org/10.4418/2020.75.1.16}{\nolinkurl{10.4418/2020.75.1.16}}}
}

@article{Simon1993,
	author = {Simon, Leon},
	title = {Existence of surfaces minimizing the {Willmore} functional},
	journal = {Comm. Anal. Geom.},
	volume = {1}, number = {2}, year = {1993}, pages = {281--326},
	doi = {10.4310/CAG.1993.v1.n2.a4},
	note = {DOI \href{https://doi.org/10.4310/CAG.1993.v1.n2.a4}{\nolinkurl{10.4310/CAG.1993.v1.n2.a4}}}
}

@article{Topping1997,
	author = {Topping, Peter M.},
	title = {The optimal constant in {Wente}'s {$L^\infty$} estimate},
	journal = {Comment. Math. Helv.},
	volume = {72}, number = {2}, year = {1997}, pages = {316--328},
	doi = {10.1007/s000140050018},
	note = {DOI \href{https://doi.org/10.1007/s000140050018}{\nolinkurl{10.1007/s000140050018}}}
}

@article{Topping2008,
	author = {Topping, Peter M.},
	title = {Relating diameter and mean curvature for submanifolds of {Euclidean} space},
	journal = {Comment. Math. Helv.},
	volume = {83}, number = {3}, year = {2008}, pages = {539--546},
	doi = {10.4171/CMH/135},
	note = {DOI \href{https://doi.org/10.4171/CMH/135}{\nolinkurl{10.4171/CMH/135}}}
}

@article{Helmers2015,
	author = {Helmers, Michael},
	title = {Convergence of an approximation for rotationally symmetric two-phase lipid bilayer membranes},
	journal = {Quart. J. Math.},
	volume = {66}, number = {1}, year = {2015}, pages = {143--170},
	doi = {10.1093/qmath/hau027},
	note = {DOI \href{https://doi.org/10.1093/qmath/hau027}{\nolinkurl{10.1093/qmath/hau027}}}
}

@article{Miura2021,
	author = {Miura, Tatsuya},
	title = {Geometric inequalities involving mean curvature for closed surfaces},
	journal = {Selecta Math. (N.S.)},
	volume = {27}, number = {5}, year = {2021}, pages = {Paper No. 80, 24 pp.},
	doi = {10.1007/s00029-021-00696-5},
	note = {DOI \href{https://doi.org/10.1007/s00029-021-00696-5}{\nolinkurl{10.1007/s00029-021-00696-5}}}
}

@article{Paeng2014,
	author = {Paeng, Seong-Hun},
	title = {Diameter of an immersed surface with boundary},
	journal = {Differential Geom. Appl.},
	volume = {33}, year = {2014}, pages = {127--138},
	doi = {10.1016/j.difgeo.2014.02.007},
	note = {DOI \href{https://doi.org/10.1016/j.difgeo.2014.02.007}{\nolinkurl{10.1016/j.difgeo.2014.02.007}}}
}

@article{MenneScharrer2026,
	author = {Menne, Ulrich and Scharrer, Christian},
	title = {A priori bounds for geodesic diameter. {Part III}. {A} {Sobolev--Poincar\'e} inequality and applications to a variety of geometric variational problems},
	journal = {Duke Math. J.},
	volume = {175}, number = {9}, year = {2026}, pages = {1649--1714},
	doi = {10.1215/00127094-2025-0051},
	note = {DOI \href{https://doi.org/10.1215/00127094-2025-0051}{\nolinkurl{10.1215/00127094-2025-0051}}}
}

@book{GolubitskyGuillemin1973,
	author = {Golubitsky, Martin and Guillemin, Victor},
	title = {Stable Mappings and Their Singularities},
	series = {Graduate Texts in Mathematics},
	volume = {14},
	publisher = {Springer-Verlag},
	address = {New York},
	year = {1973},
	doi = {10.1007/978-1-4615-7904-5},
	note = {DOI \href{https://doi.org/10.1007/978-1-4615-7904-5}{\nolinkurl{10.1007/978-1-4615-7904-5}}}
}

@article{Minkowski1903,
	author = {Minkowski, Hermann},
	title = {Volumen und {O}berfl\"ache},
	journal = {Math. Ann.},
	volume = {57},
	number = {4},
	year = {1903},
	pages = {447--495},
	doi = {10.1007/BF01445180},
	note = {DOI \href{https://doi.org/10.1007/BF01445180}{\nolinkurl{10.1007/BF01445180}}}
}

@article{WuZheng2011,
	author = {Wu, Jia-Yong and Zheng, Yu},
	title = {Relating diameter and mean curvature for {Riemannian} submanifolds},
	journal = {Proc. Amer. Math. Soc.},
	volume = {139}, number = {11}, year = {2011}, pages = {4097--4104},
	doi = {10.1090/S0002-9939-2011-10848-7},
	note = {DOI \href{https://doi.org/10.1090/S0002-9939-2011-10848-7}{\nolinkurl{10.1090/S0002-9939-2011-10848-7}}}
}

@article{FlaimScharrer2024,
	author = {Flaim, Marco and Scharrer, Christian},
	title = {Diameter estimates for surfaces in conformally flat spaces},
	journal = {Manuscripta Math.},
	volume = {174}, number = {3--4}, year = {2024}, pages = {1005--1014},
	doi = {10.1007/s00229-024-01539-1},
	note = {DOI \href{https://doi.org/10.1007/s00229-024-01539-1}{\nolinkurl{10.1007/s00229-024-01539-1}}}
}

@article{ChambersMarxKuo2026,
	author = {Chambers, Gregory R. and Marx-Kuo, Jared},
	title = {Diameter, area, and mean curvature},
	journal = {Differential Geom. Appl.},
	volume = {103}, year = {2026}, pages = {Paper No. 102376},
	doi = {10.1016/j.difgeo.2026.102376},
	note = {DOI \href{https://doi.org/10.1016/j.difgeo.2026.102376}{\nolinkurl{10.1016/j.difgeo.2026.102376}}}
}

@article{Wu2023,
	author = {Wu, Jia-Yong},
	title = {Diameter estimates for submanifolds in manifolds with nonnegative curvature},
	journal = {Differential Geom. Appl.},
	volume = {90}, year = {2023}, pages = {Paper No. 102048},
	doi = {10.1016/j.difgeo.2023.102048},
	note = {DOI \href{https://doi.org/10.1016/j.difgeo.2023.102048}{\nolinkurl{10.1016/j.difgeo.2023.102048}}}
}

@article{LiuTaiWu2025,
	author = {Liu, Leyi and Tai, Shengqiang and Wu, Jia-Yong},
	title = {Diameter bounds for submanifolds in manifolds with asymptotically nonnegative curvature},
	journal = {Results Math.},
	volume = {80}, number = {5}, year = {2025}, pages = {Paper No. 156},
	doi = {10.1007/s00025-025-02471-w},
	note = {DOI \href{https://doi.org/10.1007/s00025-025-02471-w}{\nolinkurl{10.1007/s00025-025-02471-w}}}
}
\end{document}